\documentclass{amsart}

\newtheorem{theorem}{Theorem}[section]
\newtheorem{lemma}[theorem]{Lemma}

\theoremstyle{definition}
\newtheorem{definition}[theorem]{Definition}

\newtheorem{proposition}[theorem]{Proposition}
\theoremstyle{remark}
\newtheorem{remark}[theorem]{Remark}

\numberwithin{equation}{section}

\begin{document}

\title{Non-stationary $\alpha$-fractal functions and their dimensions in various function spaces}



\author{Anarul Islam Mondal}
\address{Department of Mathematics, NIT Rourkela, India 769008}
\email{anarulmath96@gmail.com}
\author{Sangita Jha}
\address{Department of Mathematics, NIT Rourkela, India 769008}
\email{jhasa@nitrkl.ac.in}




\keywords{Fractal functions (primary) \and Attractor \and Non-stationary iterated function system, \and Function spaces \and Fractal dimension}
\subjclass{28A80 (primary) 26A18 \and 35B41 \and 41A30 \and 46B70} 
\begin{abstract}
In this article, we study the novel concept of non-stationary iterated function systems (IFSs) introduced by Massopust in 2019. At first, using a sequence of different contractive operators, we construct non-stationary $\alpha$-fractal functions on the space of all continuous functions. Next, we provide some elementary properties of the fractal operator associated with the non-stationary $\alpha$-fractal functions. Further, we show that the proposed interpolant generalizes the existing stationary interpolant in the sense of IFS. For a class of functions defined on an interval, we derive conditions on the IFS parameters so that the corresponding non-stationary $\alpha$-fractal functions are elements of some standard spaces like bounded variation space, convex Lipschitz space, and other function spaces. Finally, we discuss the dimensional analysis of the corresponding non-stationary $\alpha$-fractal functions on these spaces.
\end{abstract}

\maketitle



\section{INTRODUCTION}
Traditional interpolants, such as polynomials, trigonometric, rational, and spline functions, are always differentiable many times, with the possible exception of a finite set of points.  On the other hand, real-world problems are complicated and rarely exhibit a sense of smoothness in their traces. Fractal functions are designed to approximate complicated, extremely
irregular structures. However, despite its natural appearance, this approach received significantly less attention until 1986. Barnsley \cite{1,2} first introduced fractal interpolation functions (FIFs) to address this issue. FIFs, in general, are self-similar/affine and the Hausdorff Besicovitch dimensions of their graphs are non-integers. The main advantage here is the free choice of scaling factors and the self-referentiality character of fractal functions. The free choice allows us to select either a smooth or a non-smooth approximant. One can generalize classical interpolation techniques using smooth FIFs, for instance, see \cite{3}. 

Inspired by Barnsley's construction of FIFs and targeting the non-smooth approximants, Navascu\'es \cite{6}  introduced a family of fractal functions associated with a given continuous $f$ and the IFS parameters. This method of fractal perturbation provides a bounded linear operator, known as the $\alpha$-fractal operator \cite{Maria,12}, which links the theory of FIF to the area of Functional analysis, Operator theory, Harmonic analysis, and Approximation theory. 
Also,  many researchers have looked into the theory of dimensional and analytical aspects of FIFs in various directions and domains (for example, see the contribution of Chand\cite{3}, Viswanathan\cite{12}, Vijender\cite{Vij1}, Verma\cite{15}, Ruan\cite{Ruan1}, and others \cite{Vishal,6,Prasad,sahu} 
and the references therein ).

A useful method to construct fractals is by obtaining the fixed points of contractive operators for a special type of IFS \cite{H}. In the existing literature, the fractal defined as the attractor of a single IFS is self-similar, that is, its local shape is consistent under certain contraction maps. However, it has been observed that a sequence of IFSs is used in non-stationary subdivision schemes \cite{Dyn}. 
Recently, Levin, et al. \cite{10} have introduced a broader category of sequences comprising several contractive operators. As a generalisation of the Banach fixed point theorem, they study the trajectories of contraction mappings in \cite{Dyn,10}. It creates limit attractors with varying shapes or features at various scales. Up to now, researchers have used one contractive operator and iterated it finite or infinite times to get a stationary fractal function, which may not always provide a new class of fractal functions. By utilising the idea of forward and backward trajectories, Massopust \cite{11} introduces new types of fractal functions with various local and global behaviours  and expands fractal interpolation to a new and more adaptable environment.

In this paper, we define the aforementioned fractal operator on the set of all continuous functions but now in the non-stationary setting. As we use different contractive operators, this may give new fractal functions. We also study the essential properties of such an operator. 
The study of the fractal operators in various function spaces helps in investigating shape-preserving approximation in those function spaces. We continue to explore the aforementioned fractal operator with the non-stationary setting on the space of functions of bounded variations $\mathcal{\mathcal{BV}}(I)$, function space $V_\beta ([0,1])$ and on the convex Lipschitz space $\mathcal{V}^{\theta}(I)$. 

The study of computing the dimension of fractal sets is one of the open problems in fractal geometry. Recently, several researchers made serious efforts to compute the box/Hausdorff dimension of fractal functions and $\alpha$-fractal functions. We refer the reader to \cite{Nasim,18} for studying the fractal dimension of stationary $\alpha$-fractal functions in a few function spaces and in \cite{Ruan2} for studying the box-dimension of general recurrent fractal functions. An effort on calculating the fractal dimension of the Riemann-Liouville fractional integral of 1-dimensional continuous functions was made by Liang \cite{17,Liang}. In the present article, we attempt to find a bound of the box dimension of the proposed interpolant by constructing the non-stationary $\alpha$-fractal functions in suitable function spaces. Also, we point out that our results are generalizing certain existing results for appropriate parameters.

FIFs have been found to have greater advantages than traditional interpolants when it comes to fitting and approximating naturally occurring functions with self-similarity. In practice,
the FIF method has been used in disciplines such as image compression \cite{Ali}, signal processing \cite{David}, and physics \cite{Basu} as an alternative to traditional interpolation methods. The stationary FIF can have local or global data point dependence, and FIF maintains self-referentiality. In addition, the non-stationary settings advanced fractal functions to incorporate the scale and location dependent features also. These motivate the study of non-stationary $\alpha$-fractal functions and we believe that the work in the present article will also find many applications to the best of our knowledge.

The rest of the article is organized as follows. We review the concepts of trajectories and IFSs that are necessary to build non-stationary FIFs in Section 2. We describe the construction of the non-stationary $\alpha$-fractal functions in Section 3. The related fractal operator on $\mathcal{C}(I)$ is explored in Section 4. In the final part, we define non-stationary $\alpha$-fractal functions on various function spaces and study their dimensions.

\section{Notation and Preliminaries}
For a fixed $ k \in \mathbb{N}$, we shall write \\
$$\mathbb{N}_{k} = \{ 1, 2, 3,\dots, k \}\ \textit{and}\ \ \mathbb{N}_{k}^{0} = \{0, 1, 2, 3, \dots, k \}.$$\\
Let $(X,d)$ be a complete metric space. For a map $w: X \longrightarrow X$, let 
$$ Lip(w) = \sup \left\{ \frac{d(w(x),w(y))}{d(x,y)} : x,y \in X , x \neq y \right\}$$
denote the Lipschitz constant of $w$. If $Lip(w) < \infty $, then $w$ is called a Lipschitz function, and if $Lip(w ) < 1$, then $w$ is called a contraction. \\
Let $\mathcal{H}(X)$ denote the collection of all non-empty compact subsets of $X$. For  $C_1,C_2\in \mathcal{H}(X)$, define their Hausdorff distance as 
$$h(C_1,C_2) = max \{ d(C_1,C_2) , d(C_2,C_1) \},$$
where  $d(C_1,C_2)= \sup\limits_{x \in C_1} \inf\limits_{y \in C_2} d(x,y)$.
The space $(\mathcal{H}(X), h )$ is a complete metric space known as the space of fractals.

\begin{definition}
An iterated function system(IFS) $\mathcal{I} = \{ X ; w_{i} : i\in \mathbb{N}_N\}$ consists of a complete metric space $(X,d)$ with $N$ continuous maps $ w_{i}:X \longrightarrow X$.\\
The IFS $\mathcal{I}$ is hyperbolic if each $w_{i}$ in $\mathcal{I}$ is a contraction.
\end{definition}
For a hyperbolic IFS $\mathcal{I}$, the set valued Hutchinson map $W: \mathcal{H}(X) \longrightarrow \mathcal{H}(X)$ is defined as 
$$W(B) = \bigcup_{i=1}^{N} w_{i}(B).$$ 
It is known that $W$ is a contraction map on $\mathcal{H}(X)$ with the Lipschitz constant
$Lip(W) = \max \{ Lip(w_{i}):i=1,2,\dots,N\}.$ By using the Banach fixed point theorem, there exists a unique $A \in \mathcal{H}(X)$ such that $A = W(A)$. This $A$ is called the attractor of the IFS. The attractor $A$ can be obtained as the limit of the iterative process $A_{k} = W(A_{k-1}); k \in \mathbb{N}$, where $A_{0} \in \mathcal{H}(X)$ is any arbitrary set.
\\\\
Notice that as $A$ satisfies the self-referential equation 
$$ A = W(A) = \bigcup_{i=1}^{N} w_{i} (A) ,$$
the attractor is, in general, a fractal set.
\\\\
Let $(X,d)$ be a complete metric space and $\{T_{m}\}_{m \in \mathbb{N}}$  be a sequence of transformations on $X $.
\begin{definition}
A subset $\mathcal{P}$ of $X$ is called an invariant set of the
sequence $\{T_{m}\}_{m \in \mathbb{N}}$ if for all $ m \in \mathbb{N}$ and for all $ x \in \mathcal{P} , T_{m}(x) \in \mathcal{P}$ .
\end{definition}

We shall look at the following result to determine how to obtain an invariant set from a sequence of transformations $\{T_{m}\}_{m \in \mathbb{N}}$.

\begin{lemma}[\cite{10}]
Let $\{T_{m}\}_{m \in \mathbb{N}}$ be a sequence of transformations on $(X, d)$. Suppose there exists a $q \in X$ such that for all $x \in X$
$$d(T_{m}(x), q) \le \mu d(x, q) + M,  \mu \in [0, 1), M>0.$$ Then the ball $B_{r}(q)$ of radius $r = \frac{M}{1-\mu}$ centered at $q$ is an invariant set for $\{T_{m}\}_{m \in \mathbb{N}}$.
\end{lemma}

\begin{definition}(Forward and Backward Trajectories)
Let $\{T_{m}\}_{m \in \mathbb{N}}$ be a sequence of Lipschitz maps on the metric space $X$. The Forward and backward procedures are defined as
$$ \phi_m := T_{m}\ o\ T_{m-1}\ o\ \dots\ o\ T_{1} \ \ \text{and} \ \ \psi_m := T_{1}\ o\ T_{2}\ o\ \dots\ o\ T_{m}.$$
\end{definition}

The limits of forward trajectories might not always produce new fractal classes, as was noticed in \cite{10}. However, backward trajectories converge under relatively moderate conditions, even when forward trajectories do not converge to a (contractive) IFS, and may lead to the generation of new classes of fractal sets. We summarize the result in the following theorem.

\begin{theorem}\cite{10}
Let $\{ W_{m}\}_{m \in \mathbb{N}}$ be a family of set-valued maps of the form
\begin{equation*}
{ W_{m}(A_0) := \bigcup_{i=1}^{n_m} w_{i,m}(A_0) ,\ \ A_{0} \in \mathcal{H}(X),}
\end{equation*}
where the elements are collections $W_m = \{w_{i,m}: i \in \mathbb{N}_{n_m}\}$  of contractions constituting an IFS on the complete metric space $(X,d)$. Assume that 
\begin{enumerate}
\item there exists a nonempty closed invariant set $\mathcal{P} \subset X$ for $\{w_{i,m}\},\ i \in \mathbb{N}_{n_m} , m \in \mathbb{N};$\\
and
\item $\sum\limits_{m=1}^{\infty} \prod\limits_{j=1}^{m} Lip(W_j) < \infty.$

\end{enumerate}
Then the backward trajectories $\{\psi_{m}(A_0)\}$ converge to a unique attractor $A \subseteq \mathcal{P}$ for any initial $A_{0} \subseteq \mathcal{P}$ .
\end{theorem}

We now recall the two notions of fractal dimensions. For more details, readers are encouraged to study the book \cite{19}.
\begin{definition}
    Let $N_{\rho}(E)$ be the least number of sets with the diameter at most $\rho$ that can cover $E$, where $E$ is a non-empty bounded subset of $\mathbb{R}^n$. The upper box dimension and lower box dimension of $E$, respectively, are defined as
    $$\overline{\dim}_B (E) =\limsup_{\delta \to 0} \dfrac{\log N_{\rho}(E)}{\log \dfrac{1}{\rho}},$$
    $$\underline{\dim}_B (E) =\liminf_{\rho \to 0} \dfrac{\log N_{\rho}(E)}{\log \dfrac{1}{\rho}}.$$
\end{definition}

\begin{definition}
    The diameter of a non-empty set $U \subset \mathbb{R}^n$ is defined as
    $$|U| = \sup \{ |x-y| : x,y \in U \}.$$
    Let $E$ be a non-empty bounded subset of $\mathbb{R}^n$. We say that $\{U_i\}$ is a $\rho$-cover of $E$ if $\{U_i\}$ is a countable collection of sets of diameter at most $\rho$ which covers $E$. Let $s \ge 0$. For any $\rho > 0$, we define
    $$H_{\rho}^{s} (E) = \inf \left\{ \displaystyle \sum_{i=1}^{\infty} |U_i|^s : \{U_i\}\ \text{is a}\ \rho-\text{cover of}\ F\right\}.$$
    The $s$-dimensional Hausdorff measure of $E$ is defined by 
    $$ H^{s}(E) = \lim_{\rho \to 0} H_{\rho}^{s} (E).$$
\end{definition}

\begin{definition}
    Let $s \ge 0$. The Hausdorff dimension of a set $E \subset \mathbb{R}^n$ is defined by 
    $$\dim_H (E) = \sup \{ s : H^{s}(E)= \infty \} = \inf \{ s : H^{s}(E)= 0 \}.$$
\end{definition}
\section{Construction of non-stationary $\alpha$ -fractal function}
Let $I=[a,b]$ and $f: I \longrightarrow \mathbb{R}$ be a continuous function. Define a partition $\Delta$ by 
$$ \Delta = \{ ( x_0, x_1,\dots, x_N ) : a=x_0 < x_1 <\dots< x_N=b \}.$$
For $i \in \mathbb{N}_N$, let $I_i = [x_{i-1},x_{i}]$.
Suppose the affine maps $l_{i} : I \longrightarrow I_{i}$ are defined as follows 
$l_{i}(x) = a_{i}x+ e_{i},\;  i \in \mathbb{N}_{N},$
where $a_{i}, e_{i}$ are chosen in such a way that the maps $l_{i}$ satisfy $
    l_{i}(x_0)=x_{i-1},\; l_{i}(x_N)=x_{i}.
$ Let  $m \in \mathbb{N}$ and set $\mathbf{K} = I \times \mathbb{R}$.  We use the following notation:
\begin{center}
$\alpha_{m}:= ( \alpha_{1,m} , \alpha_{2,m},\dots,\alpha_{N,m}) ,\ \  \alpha := \{ \alpha_m \}_{m \in \mathbb{N}}$\ \ and\ \ $b:= \{ b_m \}_{m \in \mathbb{N}}.$
\end{center}

We define $F_{i,m} : \mathbf{K} \longrightarrow \mathbb{R}$ by 
\begin{equation*}
    F_{i,m}(x,y) = \alpha_{i,m}(x) y + f(l_{i}(x))-\alpha_{i,m}(x) b_m(x),
\end{equation*}
where $\alpha_{i,m} :I \longrightarrow \mathbb{R}$ are continuous functions such that $$\|\alpha\|_\infty = \sup\{\|\alpha_{m}\|_{max} : m \in \mathbb{N}\} < 1, \ \ \text{where}\ \ ||\alpha_{m}||_{max} = \max \{||\alpha_{i,m}||_\infty : i \in \mathbb{N}_N\},$$
 and $b_m \in \mathcal{C}(I)$ such that $$b_m \neq f,\  b_m(x_0)=f(x_0)\ \  and\ \ b_m(x_N)= f(x_N).$$
For each $i \in \mathbb{N}_{N}$, we define $$W_{i,m} : \mathbf{K} \longrightarrow I_{i} \times \mathbb{R} \;  \text{by} \;  W_{i,m}(x,y) = (l_{i}(x) , F_{i,m}(x,y)).$$
Now we have a sequence of IFSs
$$\mathcal{I}_m = \{ \mathbf{K};W_{i,m} : i \in \mathbb{N}_{N}\} .$$
One can show that for each $m \in \mathbb{N}$, the IFS $\mathcal{I}_m$ has a unique attractor $G$, and it is the graph of a continuous function that interpolates the given data \cite{1}.

\begin{proposition}\cite[Proposition 2.6]{13}
Let $\{T_m\}_{m \in \mathbb{N}}$ be a sequence of Lipschitz maps on a complete metric space $(X,d)$ with Lipschitz constant $\delta_m$. If there exists $ x^{*} $ in the space such that the sequence $\{ d( x^* , T_{m}(x^*))\}$ is bounded, and $\displaystyle \sum_{m=1}^{\infty} \prod_{i=1}^{m} \delta_i < \infty,$ then the sequence $\{ \psi_{m}(x) \}$ converges for all $x \in X$ to a unique limit $\bar{x}$.
\end{proposition}
Now, let $ C_{f}(I) = \{ g \in \mathcal{C}(I) : g(x_i)=f(x_i),\ i= 0,N \}$. Then $C_{f}(I)$ is a complete metric space.
For $m \in \mathbb{N}$, we define a sequence of Read-Bajraktarevi\'c (RB) operators $ T^{\alpha_m} : C_{f}(I) \longrightarrow C_{f}(I)$ by 
$$(T^{\alpha_m} g)(x) = F_{i,m}\ (\ {l_{i}}^{-1}(x) ,\ g({l_{i}}^{-1}(x))), \;  x \in I_{i},\ \ i \in \mathbb{N}_N .$$
\begin{proposition}\cite[Proposition 2.9.]{13} The above operators $ T^{\alpha_m} : C_{f}(I) \longrightarrow C_{f}(I)$ are well defined for each $m \in \mathbb{N}$.
\end{proposition}
Using similar ideas from \cite{14,13}, we have the following result.
\begin{theorem}\label{mainthm}
Consider the sequence of operators $\{ T^{\alpha_m} \}_{m \in \mathbb{N}}$ on $ C_f (I)$ defined above with the conditions described. Then for every $h \in C_f (I)$, the sequence\\
$\{ T^{\alpha_1}\ o\ T^{\alpha_2}\ o \dots o\ T^{\alpha_m} h \}$ converges to a map $f_{b}^{\alpha}$ of $C_f (I)$.
\end{theorem}
\begin{definition}
The function $f_{b}^\alpha$ is called a non-stationary $\alpha$-fractal function with respect to $f, \alpha, b$ and the partition $\Delta$ as described above.
\end{definition}
\begin{remark}
Note that, as each $ T^{\alpha_m} $ is a contraction, there is a unique stationary function $f_{m}^{\alpha}$ such that $ T^{\alpha_m}(f_{m}^{\alpha}) = f_{m}^{\alpha} $ and it satisfies the functional equation:
$$f_{m}^{\alpha}(x) = F_{i,m}\ (\ Q_{i}(x) ,\ f_{m}^{\alpha}(Q_{i}(x)))\ \ \forall\ \ x \in I_i,$$
where $Q_{i}(x) := {l_{i}}^{-1}(x)$. That is,
$$ f_{m}^{\alpha}(x) = f(x) + \alpha_{i,m}(Q_{i}(x)). f_{m}^{\alpha}(Q_{i}(x))- \alpha_{i,m}(Q_{i}(x)) b_{m}(Q_{i}(x)).$$
\end{remark}

\section{Associated fractal operator on $\mathcal{C}(I)$}
Let $||\alpha||_\infty =\sup\limits_{m \in \mathbb{N}} ||\alpha_m||_{max} < 1$ and $||b||_\infty: =\sup\limits_{m \in \mathbb{N}} ||b_m||_\infty<\infty.$
We consider $b_m = L_{m}f$ such that  $L_m : \mathcal{C}(I) \to \mathcal{C}(I)$ is a linear bounded operator satisfying $L_{m} f(x_i) = f(x_i)$ for $m \in \mathbb{N},\ i=0,N,$ and $||L||_\infty:=\sup\limits_{m \in \mathbb{N}} ||L_m|| < \infty.$ Let $f \in \mathcal{C}(I)$. We define the $\alpha$-fractal operator $\mathcal{F}_{b}^{\alpha} \equiv \mathcal{F}_{\Delta,b}^{\alpha}$ as 
\begin{center}
$\mathcal{F}_{b}^{\alpha} : \mathcal{C}(I) \to \mathcal{C}(I) ,\ \ \mathcal{F}_{b}^{\alpha}(f) = f_{b}^{\alpha}.$
\end{center}

\begin{lemma}\label{lemma}
Let $X$ be a Banach space and $T: X \to X$ be a bounded linear operator. If $||T|| < 1$, $(Id - T)^{-1}$ exists and bounded, where $Id$ denotes the identity operator on $X$.
\end{lemma}

\begin{lemma}\label{lemma3}
Let $X$ be a normed linear space and $T: X \to X$ be a bounded linear operator, and  $S: X\to X$ be a compact operator. Then $ST$ and $TS$ are compact operators.
\end{lemma}

We now describe a few properties of the non-stationary $\alpha$-fractal operator. Note that the following properties are also studied in the literature for the stationary case \cite{6,20}. However, the approach in the non-stationary setting is different, and for better understanding and completeness, we study the following properties.
Note that we have considered $L_m\in \mathcal{C}(I)$ and $||L||_\infty:=\sup\limits_{m \in \mathbb{N}} ||L_m|| < \infty $. Hence $C_L:=\sup\limits_{m \in \mathbb{N}} \{||Id-L_m||\}<\infty$.
\begin{theorem}
Let $||\alpha||_\infty <1$ and $Id$ be the identity operator on $\mathcal{C}(I)$.
\begin{enumerate}
\item For $f \in \mathcal{C}(I)$, the perturbation error satisfies the following inequality:
$$||f_{b}^{\alpha} - f||_\infty \le \frac{||\alpha||_\infty}{1-||\alpha||_\infty} \sup_{m \in \mathbb{N}}\{ ||f-L_{m}f||_\infty\} \le \frac{||\alpha||_\infty}{1-||\alpha||_\infty} C_L ||f||_\infty.$$
\item If $\alpha=0$, then $\mathcal{F}_{b}^{\alpha}$ is norm preserving. Infact it holds that $\mathcal{F}_{b}^{0} \equiv Id$.
\item The fractal operator $ \mathcal{F}_{b}^{\alpha}: \mathcal{C}(I) \longrightarrow \mathcal{C}(I)$ is linear and bounded with respect to the uniform norm.
\item For a suitable value of the scaling function, the operator $\mathcal{F}_{b}^{\alpha}$ is an approximation type operator.
\item For $||\alpha||_\infty < \frac{1}{1+C_L}$ , the fractal operator $\mathcal{F}_{b}^{\alpha}$ is bounded below. In particular, $\mathcal{F}_{b}^{\alpha}$ is one to one.
\item If $||\alpha||_\infty < \frac{1}{1+C_L}$, then $\mathcal{F}_{b}^{\alpha}$ has a bounded inverse and consequently a topological isomorphism.
\item The fixed points of $L_m$ are also the fixed points of $\mathcal{F}_{b}^{\alpha}$.
\item If $1$ belongs to the spectrum of $L_m$, then $ 1 \le ||\mathcal{F}_{b}^{\alpha}||$.
\item For $||\alpha||_\infty < \frac{1}{1+C_L}$ , the fractal operator $\mathcal{F}_{b}^{\alpha}$ is not a compact operator.
\item For $||\alpha||_\infty < \frac{1}{1+C_L}$ , the fractal operator $\mathcal{F}_{b}^{\alpha}$ has a closed range.
\end{enumerate}
\end{theorem}
\begin{proof}
\begin{enumerate}
\item From the definition of RB operators, we have 
\begin{equation}\label{eq-RB}
(T^{\alpha_m}g)(x) = f(x) + \alpha_{i,m}(Q_{i}(x)). g(Q_{i}(x))- \alpha_{i,m}(Q_{i}(x)) b_{m}(Q_{i}(x))
\end{equation}
for all $ x \in I_i,\ \ i \in \mathbb{N}_N\ \text{and}\ m \in \mathbb{N}.$\\
Now $ \forall\ x \in I_i,$
$$ T^{\alpha_1}\ o\ T^{\alpha_2}\ o\dots o\ T^{\alpha_m} f(x) - f(x) = \alpha_{i,1}(Q_{i}(x)) ( T^{\alpha_2}\ o\ T^{\alpha_3}\ o \dots o\ T^{\alpha_m} f - b_1)(Q_{i}(x)).$$
Inductively, we get $\forall\ x \in I_i$,
\begin{equation}\label{eq-operator}
T^{\alpha_1}\ o\ T^{\alpha_2}\ o \dots o\ T^{\alpha_m} f(x) - f(x) = \sum_{l=1}^{m} \alpha_{i,1}(Q_{i}(x))\dots\alpha_{i,l}(Q_{i}^{l}(x)) (f-b_l)(Q_{i}^{l}(x),
\end{equation}
where $Q_{i}^{l}$ is a suitable finite composition of maps $Q_i$.\\
Taking limit as $m \rightarrow \infty$, we get
\begin{align}\label{eq-functional}
f_{b}^{\alpha} (x) - f(x) &= \lim_{m \rightarrow \infty} \sum_{l=1}^{m} \alpha_{i,1}(Q_{i}(x))\dots\alpha_{i,l}(Q_{i}^{l}(x)) (f-b_l)(Q_{i}^{l}(x))\nonumber \\
&= \lim_{m \rightarrow \infty} \sum_{l=1}^{m} \alpha_{i,1}(Q_{i}(x))\dots\alpha_{i,l}(Q_{i}^{l}(x)) (f-L_l f)(Q_{i}^{l}(x)).
\end{align}
So,
\begin{equation*}
\begin{split}
||f_{b}^{\alpha}- f||_\infty & \le \lim_{m \rightarrow \infty} \sum_{l=1}^{m} ||\alpha||_{\infty}^{l} ||f-L_l f||_\infty \\
& \le \lim_{m \rightarrow \infty} \sum_{l=1}^{m} ||\alpha||_{\infty}^{l} \sup_{m \in \mathbb{N}}||f-L_m f||_\infty\\
&= \sum_{l=1}^{\infty} ||\alpha||_{\infty}^{l} \sup_{m \in \mathbb{N}}||f-L_{m}f||_\infty\\
& = \frac{||\alpha||_\infty}{1-||\alpha||_\infty}\sup_{m \in \mathbb{N}}||f-L_m f||_\infty\\
& = \frac{||\alpha||_\infty}{1-||\alpha||_\infty} C_L ||f||_\infty.
\end{split}
\end{equation*}
\begin{equation}\label{eqbound}
||f_{b}^{\alpha}- f||_\infty \le \frac{||\alpha||_\infty}{1-||\alpha||_\infty} C_L ||f||_\infty.
\end{equation}

\item From equation \eqref{eqbound}, we have
$$||f_{b}^{\alpha} - f||_\infty \le \frac{||\alpha||_\infty}{1-||\alpha||_\infty} C_L ||f||_\infty$$
$$\implies ||\mathcal{F}_{b}^{\alpha}(f) - f||_\infty \le \frac{||\alpha||_\infty}{1-||\alpha||_\infty} C_L ||f||_\infty.$$
If $\alpha = 0 $, then $||\mathcal{F}_{b}^{\alpha}(f) - f||_\infty = 0 $. Therefore,  $\mathcal{F}_{b}^{\alpha}(f) = f\ \ \implies \mathcal{F}_{b}^{0} = Id.$

\item Let $f,g \in \mathcal{C}(I)$ and $c,d \in \mathbb{R}$. Then from equation \eqref{eq-functional}, we have for all $x \in I_i$, 
$$(cf)_{b}^{\alpha} (x) = (cf)(x) + \lim_{m \rightarrow \infty} \sum_{l=1}^{m} \alpha_{i,1}(Q_{i}(x))\dots\alpha_{i,l}(Q_{i}^{l}(x)) (cf-L_l (cf))(Q_{i}^{l}(x)),$$
$$(dg)_{b}^{\alpha} (x) = (dg)(x) + \lim_{m \rightarrow \infty} \sum_{l=1}^{m} \alpha_{i,1}(Q_{i}(x))\dots\alpha_{i,l}(Q_{i}^{l}(x)) (dg-L_l (dg))(Q_{i}^{l}(x)).$$
As $L_l$ is linear, so that 
\begin{align*}
    &(cf)_{b}^{\alpha} (x) + (dg)_{b}^{\alpha} (x)\\ &= (cf+dg)(x) + \lim_{m \rightarrow \infty} \sum_{l=1}^{m} \alpha_{i,1}(Q_{i}(x))\dots\alpha_{i,l}(Q_{i}^{l}(x)) [cf+dg-L_l (cf+dg)](Q_{i}^{l}(x)).
\end{align*}

Also,
\begin{align*}
(cf+dg)_{b}^{\alpha} (x) &= (cf+dg)(x)\\ &+ \lim_{m \rightarrow \infty} \sum_{l=1}^{m} \alpha_{i,1}(Q_{i}(x))\dots\alpha_{i,l}(Q_{i}^{l}(x)) \times
[cf+dg-L_l (cf+dg)](Q_{i}^{l}(x)).
\end{align*}
Hence we deduce that,
$$(cf+dg)_{b}^{\alpha} (x) = (cf)_{b}^{\alpha} (x) + (dg)_{b}^{\alpha} (x).$$
That is,\ \ $(cf+dg)_{b}^{\alpha} = (cf)_{b}^{\alpha} + (dg)_{b}^{\alpha} \ \ \implies \mathcal{F}_{b}^{\alpha}(cf+dg) =c \mathcal{F}_{b}^{\alpha}(f) + d \mathcal{F}_{b}^{\alpha}(g).$
This proves the linearity of the operator $\mathcal{F}_{b}^{\alpha}$.

From \eqref{eqbound}, we have
$$||f_{b}^{\alpha}||_\infty -||f||_\infty \le ||f_{b}^{\alpha}- f||_\infty \le \dfrac{||\alpha||_\infty}{1-||\alpha||_\infty} C_L ||f||_\infty.$$
That is, $$||\mathcal{F}_{b}^{\alpha}(f)||_\infty \le \left( 1 + \frac{||\alpha||_\infty}{1-||\alpha||_\infty} C_L \right)||f||_\infty.$$
$$\implies ||\mathcal{F}_{b}^{\alpha}||\le \left( 1 + \frac{||\alpha||_\infty}{1-||\alpha||_\infty} C_L \right).$$
Therefore the operator $\mathcal{F}_{b}^{\alpha}$ is bounded.

\item Let $\epsilon>0$. We choose the scaling sequence $\alpha$ such that $||\alpha||_\infty < \dfrac{\epsilon}{\epsilon + C_L ||f||_\infty}$. Using equation \eqref{eqbound}, we obtain
$$||f_{b}^{\alpha} - f ||_\infty < \epsilon$$
$$\implies ||\mathcal{F}_{b}^{\alpha}(f) - f ||_\infty < \epsilon.$$
Consequently, the operator $\mathcal{F}_{b}^{\alpha}$ is of approximation type.

\item If $||\alpha||_\infty < \frac{1}{1+C_L}$, then we have from equation \eqref{eqbound}
$$||f||_\infty - ||f_{b}^{\alpha}||_\infty \le ||f_{b}^{\alpha}- f||_\infty \le \frac{||\alpha||_\infty}{1-||\alpha||_\infty} C_L ||f||_\infty.$$
$$\implies \left( 1 - \frac{||\alpha||_\infty}{1-||\alpha||_\infty} C_L \right) ||f||_\infty \le ||\mathcal{F}_{b}^{\alpha}(f)||_\infty.$$
$$\implies \left( \frac{1 - ||\alpha||_\infty (1+C_L)}{1-||\alpha||_\infty} \right) ||f||_\infty \le ||\mathcal{F}_{b}^{\alpha}(f)||_\infty.$$
This shows that $\mathcal{F}_{b}^{\alpha}$ is bounded from below. Consequently, $\mathcal{F}_{b}^{\alpha}$ is injection.

\item From equation \eqref{eqbound}, we have
$$||Id - \mathcal{F}_{b}^{\alpha}|| \le \dfrac{||\alpha||_\infty}{1-||\alpha||_\infty} C_L.$$
As $||\alpha||_\infty < \dfrac{1}{1+C_L}$, we get $||Id - \mathcal{F}_{b}^{\alpha}|| < 1 $. Since $\mathcal{F}_{b}^{\alpha}$ is bounded, $(Id - \mathcal{F}_{b}^{\alpha})$ is also bounded. Hence by using Lemma \ref{lemma}, $\mathcal{F}_{b}^{\alpha} = Id - (Id - \mathcal{F}_{b}^{\alpha}) $ is invertible and the inverse is bounded.
\item Let for each $m \in \mathbb{N}$, $f$ is fixed point of $L_m$. Then $L_m(f)=f$ for each $m \in \mathbb{N}$. From equation \eqref{eq-functional}, we have $f_{b}^{\alpha}(x) - f(x) = 0$. So that $\mathcal{F}_{b}^{\alpha}f(x) = f(x)$. This implies that $\mathcal{F}_{b}^{\alpha}(f)=f$.

\item Let $g \in \mathcal{C}(I)$ with $||g||_\infty = 1$ and $L_m f = f$. Then item (7) gives, $\mathcal{F}_{b}^{\alpha} g = g $. Consequently, $||\mathcal{F}_{b}^{\alpha} g ||_\infty = ||g||_\infty.$ From the definition of operator norm, we have $ 1 \le ||\mathcal{F}_{b}^{\alpha}||$.

\item From item (6), we have for $||\alpha||_\infty < \dfrac{1}{1+C_L}$, the operator $\mathcal{F}_{b}^{\alpha} : \mathcal{C}(I) \longrightarrow \mathcal{C}(I)$ is one-one. Note that the range space $\mathcal{F}_{b}^{\alpha}(\mathcal{C}(I))$ is infinite dimensional. We define the inverse map $(\mathcal{F}_{b}^{\alpha})^{-1} : \mathcal{F}_{b}^{\alpha}(\mathcal{C}(I)) \longrightarrow \mathcal{C}(I)$. 
For the choice of $\alpha$, we know from item (5) that $\mathcal{F}_{b}^{\alpha}$ is bounded below and hence it follows that $(\mathcal{F}_{b}^{\alpha})^{-1}$ is a bounded linear operator. If possible, let $\mathcal{F}_{b}^{\alpha}$ be a compact operator. Then by Lemma \ref{lemma3}, we conclude that the operator $Id = (\mathcal{F}_{b}^{\alpha})(\mathcal{F}_{b}^{\alpha})^{-1}: \mathcal{F}_{b}^{\alpha}(\mathcal{C}(I)) \longrightarrow \mathcal{C}(I)$ is a compact operator, which is a contradiction to the infinite dimensionality of the space $\mathcal{F}_{b}^{\alpha}(\mathcal{C}(I))$. Hence $\mathcal{F}_{b}^{\alpha}$ is not a compact operator.

\item Let $\{ f_{b,r}^{\alpha}\}_{r \in \mathbb{N}}$ be a sequence in $\mathcal{F}_{b}^{\alpha}(\mathcal{C}(I))$ such that $f_{b,r}^{\alpha} \rightarrow g$. Then  $\{ f_{b,r}^{\alpha}\}_{r \in \mathbb{N}}$ is a Cauchy sequence in $\mathcal{C}(I)$. Thus from
$$||f_r - f_s||_\infty \le \dfrac{1 - ||\alpha||_\infty}{1-(1+C_L)||\alpha||_\infty} ||f_{b,r}^{\alpha} - f_{b,s}^{\alpha}||_\infty,$$
the sequence $\{ f_r\}$ is a Cauchy sequence in $\mathcal{C}(I)$. Consequently, there exists $f \in \mathcal{C}(I)$ such that $f_r \rightarrow f$. Using the continuity of the operator $\mathcal{F}_{b}^{\alpha}$, we obtain $ g = \mathcal{F}_{b}^{\alpha}(f_r) = f_{b,r}^{\alpha}$.
\end{enumerate}
\end{proof}
\begin{remark}
    Let $\alpha_ {i,m}(x) = \alpha_i$ and $b_m = L_m f = Lf$ for all $m \in \mathbb{N} $. Then from equation \eqref{eq-operator}, we have 
    \begin{equation*}
(T^{\alpha})^m f(x) - f(x) = \sum_{l=1}^{m} \alpha_{i}^l (f-Lf)(Q_{i}^{l}(x)).
\end{equation*}
Taking limit as $m \to \infty$, we get 
\begin{equation*}
f^{\alpha} (x) - f(x) = \lim_{m \rightarrow \infty} \sum_{l=1}^{m} \alpha_{i}^l (f-Lf)(Q_{i}^{l}(x)).
\end{equation*}
So,
\begin{equation*}
\begin{split}
||f^{\alpha}- f||_\infty & \le \lim_{m \rightarrow \infty} \sum_{l=1}^{m} |\alpha|_{\infty}^{l} ||f-L f||_\infty \\
& = (\sum_{l=1}^{\infty} |\alpha|_{\infty}^{l}) \cdot ||f-L f||_\infty\\
& \le \frac{|\alpha|_\infty}{1-|\alpha|_\infty} (1+||L||) ||f||_\infty\\
\end{split}
\end{equation*}
$\implies ||f^{\alpha}||_\infty - ||f||_\infty\ \le\ ||f^{\alpha}- f||_\infty\  \le\ \dfrac{|\alpha|_\infty}{1-|\alpha|_\infty} (1+||L||) ||f||_\infty$\\
So that, 
$$||f^{\alpha}||_\infty \le \dfrac{1 + |\alpha|_\infty ||L||}{1 - |\alpha|_\infty} \cdot ||f||_\infty,$$
which is the same upper bound of the stationary $\alpha$-fractal functions given in Proposition 2 of \cite{20}.
\end{remark}

\section{Non-stationary fractal functions on different function spaces}
The aim of this section is to study the non-stationary $\alpha$-fractal functions in different function spaces. We start this with the bounded variation space.
\subsection{Space of functions of bounded variation on $I$}
\begin{definition}
Let $f: I \to \mathbb{R}$ be a function. For each partition $P_I: x_0 < x_1 < \dots < x_N$ of the interval $I$, we define
 $$V_{P_I}(f,I) = \sum_{i=1}^{N} |f(x_i) - f(x_{i-1})|$$
 and
 $$V(f,I) = \sup_{P_I} V_{P_I}(f,I) = \sup_{P_I} \sum_{i=1}^{N} |f(x_i) - f(x_{i-1})|,$$
 where the supremum is taken over all the partitions $P_I$ of the interval $I$. If the total variation of $f$ is finite, i.e., $V(f, I) < \infty$,  we say $f$ is of bounded variation on $I$.
\end{definition} 
\begin{theorem}\cite{17}\label{thm-boundedvariation}
    If $f$ is continuous bounded variation function on an interval $I$, then $dim_H (Graph\ (f)) = dim_B (Graph\ (f)) = 1.$
\end{theorem}
Let $\mathcal{BV}(I)$ denote the set of all functions of bounded variation on $I$. On $\mathcal{BV}(I)$, we define a norm $||.||_{\mathcal{BV}}$  by
 $$||f||_{\mathcal{BV}} = |f(x_0)| + V(f,I).$$
 We know that with respect to this norm, $\mathcal{BV}(I)$ is a Banach space. Let $f \in \mathcal{BV}(I)$ and define
$\mathcal{BV}_f(I) = \{ g \in \mathcal{BV}(I) :\ f(x_0) = g(x_0),\ f(x_N) = g(x_N) \} .$ Note that $\mathcal{BV}_f(I)$ is also complete with respect to the metric induced by the norm $||.||_{\mathcal{BV}}$.
Let
$\alpha_{m}:= ( \alpha_{1,m} , \alpha_{2,m},\dots,\alpha_{N,m}) ,\ \  \alpha := \{ \alpha_m \}_{m \in \mathbb{N}},\; \; b:= \{ b_m \}_{m \in \mathbb{N}},$ and $||\alpha||_\infty =\sup\limits_{m \in \mathbb{N}} ||\alpha_m||_{max} < 1,\; \; ||b||_{\mathcal{BV}}:=\sup\limits_{m \in \mathbb{N}} ||b_m||_{\mathcal{BV}}<\infty.$ 
With all these setups, we have the following result:
\begin{theorem}\label{th-BV}
Let $ f,b_m \in \mathcal{BV}(I)$ be such that $ b_m(x_0)=f(x_0),\; b_m(x_N)=f(x_N)$ and the sequence of scaling functions $\alpha_{i,m} \in \mathcal{BV}(I)$ be such that $||\alpha_{i,m}||_\mathcal{BV} < \dfrac{1}{2N}$. Then the following hold.
\begin{enumerate}
\item The RB operator $T^{\alpha_m}$ defined in equation \eqref{eq-RB} is well defined on $\mathcal{BV}_f(I)$.
\item $T^{\alpha_m} : \mathcal{BV}_f(I) \longrightarrow \mathcal{BV}_f (I) \subset \mathcal{BV}(I)$ is, in reality, a contraction map.
\item There exists a unique function $f_{b, \mathcal{BV}}^{\alpha} \in \mathcal{BV}_f(I)$ such that the backward trajectories $ T^{\alpha_1}\ o\ T^{\alpha_2}\ o \dots o\ T^{\alpha_m} g$ of $ (T^{\alpha_m})$ converges to the map $f_{b, \mathcal{BV}}^{\alpha}$ for every $g \in \mathcal{BV}_f(I)$.
\end{enumerate}
Furthermore, we have $dim_H (Graph\ (f_{b}^{\alpha})) = dim_B (Graph\ (f_{b}^{\alpha})) = 1.$ 
\end{theorem}
\begin{proof}
\begin{enumerate}
\item From the definition of RB operators, we have for $m\in \mathbb{N}$
$$ (T^{\alpha_m}g)(x) = f(x) + \alpha_{i,m}(Q_{i}(x)). g(Q_{i}(x))- \alpha_{i,m}(Q_{i}(x)) b_{m}(Q_{i}(x)),\ \ \ x \in I_i,\ \ i \in \mathbb{N}_N.$$
As $f,\alpha_{i,m},g, b_m \in \mathcal{BV}_f(I)$, so that $ (T^{\alpha_m}g)(x) \in \mathcal{BV}_f(I)$ whenever $g \in \mathcal{BV}_f(I)$.
Therefore the RB operator $T^{\alpha_m}$ is well defined on $\mathcal{BV}_f (I)$.
\item Let $P : x_{i-1} = t_{0,i} < t_{1,i} < t_{2,i} <\dots< t_{k_i,i} = x_{i}$ be a partition of the interval $I_{i} = [x_{i-1}, x_{i}],$ where $k_i \in \mathbb{N}$. Now,
\begin{align*}
&|(T^{\alpha_m} g - T^{\alpha_m} h)(t_{j,i}) - (T^{\alpha_m} g - T^{\alpha_m} h)(t_{j-1,i})|\\
=& |\alpha_{i,m}(Q_{i}(t_{j,i}))(g-h)(Q_{i}(t_{j,i})) - \alpha_{i,m}(Q_{i}(t_{j-1,i}))(g-h)(Q_{i}(t_{j-1,i}))|\\
 \le& \bigg|\alpha_{i,m}(Q_{i}(t_{j,i}))\bigg((g-h)(Q_{i}(t_{j,i})) - (g-h)(Q_{i}(t_{j-1,i})) \bigg)\bigg|\\
&+ \bigg| \bigg(\alpha_{i,m}(Q_{i}(t_{j,i})) - \alpha_{i,m}(Q_{i}(t_{j-1,i}))\bigg)(g-h)(Q_{i}(t_{j-1,i}))\bigg|\\
\le& \|\alpha_{i,m}\|_\mathcal{BV} |(g-h)(Q_{i}(t_{j,i})) - (g-h)(Q_{i}(t_{j-1,i}))|\\
&+ |\alpha_{i,m}(Q_{i}(t_{j,i})) - \alpha_{i,m}(Q_{i}(t_{j-1,i}))| \|g-h\|_\mathcal{BV}
\end{align*}
Taking sum over $j=1$ to $ k_i$, we have
\begin{align*}
&\sum_{j=1}^{k_i} |(T^{\alpha_m} g - T^{\alpha_m} h)(t_{j,i}) - (T^{\alpha_m} g - T^{\alpha_m} h)(t_{j-1,i})|\\
\le& ||\alpha_{i,m}||_\mathcal{BV} \sum_{j=1}^{k_i} |(g-h)(Q_{i}(t_{j,i})) - (g-h)(Q_{i}(t_{j-1,i}))|\\
&+ \|g-h\|_\mathcal{BV} \sum_{j=1}^{k_i} |\alpha_{i,m}(Q_{i}(t_{j,i})) - \alpha_{i,m}(Q_{i}(t_{j-1,i}))| 
\end{align*}
Since $x_0 = Q_{i}(t_{0,i})<Q_{i}(t_{1,i})< \dots<Q_{i} (t_{k_i,i})=x_{N}$ is a partition of the interval $I=[x_0, x_N]$, so that 
\begin{align*}
    &\sum_{j=1}^{k_i} |(T^{\alpha_m} g - T^{\alpha_m} h)(t_{j,i}) - (T^{\alpha_m} g - T^{\alpha_m} h)(t_{j-1,i})|\\
    \le& ||\alpha_{i,m}||_\mathcal{BV} V(g-h, I) + \|g-h\|_\mathcal{BV} V( \alpha_{i,m} , I)\\
    \le& ||\alpha_{i,m}||_\mathcal{BV} \|g-h\|_\mathcal{BV} + \|g-h\|_\mathcal{BV} ||\alpha_{i,m}||_\mathcal{BV}\\
    =& 2 ||\alpha_{i,m}||_\mathcal{BV} \|g-h\|_\mathcal{BV}
\end{align*}

This inequality holds for every partition $P$ of $I_i$. Hence
$$ V( T^{\alpha_m} g - T^{\alpha_m} h , I_i ) \le 2 ||\alpha_{i,m}||_\mathcal{BV} \|g-h\|_\mathcal{BV}.$$
As
$$ V( T^{\alpha_m} g - T^{\alpha_m} h , I ) = \sum_{i=1}^{N} V( T^{\alpha_m} g - T^{\alpha_m} h , I_i ) \le 2N ||\alpha_{i,m}||_\mathcal{BV} \|g-h\|_\mathcal{BV}$$
and $T^{\alpha_m}g(x_0) = T^{\alpha_m}h(x_0) = f(x_0)$, we have
$$||T^{\alpha_m} g - T^{\alpha_m} h||_{\mathcal{BV}} \le 2N ||\alpha_{i,m}||_\mathcal{BV} \|g-h\|_\mathcal{BV}.$$
Since $||\alpha_{i,m}||_\mathcal{BV} < \dfrac{1}{2N}$, each $T^{\alpha_m}$ is a contraction map on the complete metric space $\mathcal{BV}_f(I)$.
\item Let $g \in \mathcal{BV}_f(I)$ be arbitrary. We check that $\{ ||T^{\alpha_m} g -g||_{\mathcal{BV}}\}$ is bounded.
To do this we first calculate $V( T^{\alpha_m} g - g , I )$. Now
\begin{align*}
&|(T^{\alpha_m} g - g)(t_{j,i}) - (T^{\alpha_m} g - g)(t_{j-1,i})|\\ 
&= |(f-g)(t_{j,i}) - (f-g)(t_{j-1,i})\\
& \hspace{0.5cm} + \alpha_{i,m}(Q_{i}(t_{j,i}))(g-b_m)(Q_{i}(t_{j,i})) - \alpha_{i,m}(Q_{i}(t_{j-1,i}))(g-b_m)(Q_{i}(t_{j-1,i}))|\\
&\le \bigg|(f-g)(t_{j,i}) - (f-g)(t_{j-1,i})\bigg|\\
& \hspace{0.5cm} + \bigg|\alpha_{i,m}(Q_{i}(t_{j,i}))\bigg((g-b_m)(Q_{i}(t_{j,i})) - (g-b_m)(Q_{i}(t_{j-1,i})) \bigg)\bigg|\\
& \hspace{0.5cm} + \bigg| \bigg(\alpha_{i,m}(Q_{i}(t_{j,i})) - \alpha_{i,m}(Q_{i}(t_{j-1,i}))\bigg)(g-b_m)(Q_{i}(t_{j-1,i}))\bigg|\\
&\le \big|(f-g)(t_{j,i}) - (f-g)(t_{j-1,i})\big| + \|\alpha_{i,m}\|_\mathcal{BV} \big|(g-b_m)(Q_{i}(t_{j,i})) - (g-b_m)(Q_{i}(t_{j-1,i}))\big|\\
& \hspace{0.5cm} + \big|\alpha_{i,m}(Q_{i}(t_{j,i})) - \alpha_{i,m}(Q_{i}(t_{j-1,i}))\big| \|g-h\|_\mathcal{BV}
\end{align*}
Taking sum over $j=1$ to $ k_i$, we have
\begin{align*}
&\sum_{j=1}^{k_i} |(T^{\alpha_m} g - g)(t_{j,i}) - (T^{\alpha_m} g - g)(t_{j-1,i})|\\
&\le \sum_{j=1}^{k_i} |(f-g)(t_{j,i}) - (f-g)(t_{j-1,i})|\\
& \hspace{0.5cm}+ ||\alpha_{i,m}||_\mathcal{BV} \sum_{j=1}^{k_i} |(g-b_m)(Q_{i}(t_{j,i})) - (g-b_m)(Q_{i}(t_{j-1,i}))|\\
& \hspace{0.5cm} + \|g-b_m\|_\mathcal{BV} \sum_{j=1}^{k_i} \big|\alpha_{i,m}(Q_{i}(t_{j,i})) - \alpha_{i,m}(Q_{i}(t_{j-1,i}))\big|\\
&\le V(f-g, I_i) + ||\alpha_{i,m}||_\mathcal{BV} V(g-b_m, I) + \|g-b_m\|_\mathcal{BV} V(\alpha_{i,m}, I)\\
&\le V(f-g, I_i) + ||\alpha_{i,m}||_\mathcal{BV} ||g-b_m||_\mathcal{BV} + \|g-b_m\|_\mathcal{BV} ||\alpha_{i,m}||_\mathcal{BV}.
\end{align*}
By a similar argument as in item (2), we get
\begin{align*}
    &V( T^{\alpha_m} g - g , I_i) \le V(f-g, I_i) + 2||\alpha_{i,m}||_\mathcal{BV} ||g-b_m||_\mathcal{BV}.
\end{align*}
Using the conditions on sequence of scaling functions, we get
\begin{align*}
    V( T^{\alpha_m} g - g , I) &= \sum_{i=1}^{N} V( T^{\alpha_m} g - g , I_i )\\
    &\le V(f-g, I) + 2 ||g-b_m||_\mathcal{BV} \sum_{i=1}^{N} ||\alpha_{i,m}||_\mathcal{BV}\\
    &\le \|f-g\|_\mathcal{BV} + \|g-b_m\|_\mathcal{BV}
\end{align*}
Also, \begin{align*}
|(T^{\alpha_m} g - g)(x_0)| &= |f(x_0) + \alpha_{1,m} (g-b_m)(x_0)-g(x_0)| = 0.
\end{align*}
Now,
\begin{align*}
&||T^{\alpha_m} g - g ||_{\mathcal{BV}}\\
&= |(T^{\alpha_m} g - g)(x_0)| + V( T^{\alpha_m} g - g , I)\\
&\le \|f-g\|_\mathcal{BV} + \|g-b_m\|_\mathcal{BV}\\
&= ||f||_{\mathcal{BV}}+||g||_{\mathcal{BV}} + ||g||_{\mathcal{BV}} + ||b_m||_{\mathcal{BV}}\\
&\le ||f||_{\mathcal{BV}}+2 ||g||_{\mathcal{BV}} + ||b||_{\mathcal{BV}}.
\end{align*}
Clearly, the bound is independent of $m$. Applying Theorem \ref{mainthm}, $\exists$ a unique $f_{b,\mathcal{BV}}^{\alpha} \in \mathcal{BV}_f(I)$ such that $f_{b,\mathcal{BV}}^{\alpha} = \displaystyle \lim_{ m \rightarrow \infty} T^{\alpha_1}\ o\ T^{\alpha_2}\ o \dots o\ T^{\alpha_m} g$ for any $g \in \mathcal{BV}_f(I)$.
Also, applying Theorem \ref{thm-boundedvariation}, we have 
$dim_H (Graph\ (f_{b}^{\alpha})) = dim_B (Graph\ (f_{b}^{\alpha})) = 1.$
\end{enumerate}
\end{proof}
\begin{remark}
    If we take $\alpha_{i,m} = \alpha_i,\; b_m=b$ for all $m \in \mathbb{N}$, then we get 
    $$f_{b,\mathcal{BV}}^{\alpha} = \displaystyle \lim_{ m \rightarrow \infty} T^{\alpha}\ o\ T^{\alpha}\ o \dots o\ T^{\alpha} g = \displaystyle \lim_{ m \rightarrow \infty} (T^{\alpha})^m$$
    for any $g \in \mathcal{BV}_f(I).$
    So that, $f_{b,\mathcal{BV}}^{\alpha} \longrightarrow f_{\mathcal{BV}}^{\alpha},\ \text{as}\ m \to \infty$, where $f_{\mathcal{BV}}^{\alpha}$ is the stationary $\alpha$-fractal function on the space $\mathcal{BV}(I)$ that appeared in \cite{Verma}. That is, the non-stationary $\alpha$-fractal function tends to stationary $\alpha$-fractal function on $\mathcal{BV}(I)$ for the above particular choice of IFS parameters.
\end{remark}
\subsection{Function space $V_{\beta}[0,1]$}
Let $\beta \in [1,2]$, and define
$$ C_\beta [0,1] =\{ f \in \mathcal{C}[0,1] : \overline{dim}_B\ G_f \le \beta \}.$$
Let $f \in \mathcal{C}[0,1]$ and $S \subset [0,1]$. We define the range of $f$ on $S$ as
$$R_f(S) = \displaystyle \sup_{x,y \in S} |f(x) - f(y)|.$$
Let $R(\delta,f) =\displaystyle \sum_{S \in \Delta_{\delta}} R_f(S)$.
For $\delta > 0$, let $\Delta_{\delta}$ be the set defined below
$$\Delta_{\delta} = \bigcup_{n=0}^{\lceil \delta^{-1} \rceil - 1}   [n \delta , (n+1)\delta].$$
It follows that \[\delta^{-1}\sum\limits_{S \in \Delta_{\delta}} R_{\lambda f}(S)\leq N_{\delta}(G_f)\leq 2(\delta^{-1}+1)+\delta^{-1}\sum\limits_{S \in \Delta_{\delta}} R_{\lambda f}(S).\]
\begin{remark} Note that the above estimation is a useful technique in fractal geometry for calculating the box dimension of fractal functions.\cite{19}
\end{remark}

\begin{theorem}
    Let $g, f \in \mathcal{C}[0,1]$ and $\lambda$ be a real number. Then for $0<\delta \le 1,$
    \begin{enumerate}
        \item[(i)] $R(\delta,\lambda f) = |\lambda| \cdot R(\delta,f)$,
        \item[(ii)] $R(\delta,f+g) \le R(\delta,f) + R(\delta,g)$,
        \item[(iii)] $R(\delta,fg) \le ||g||_\infty \cdot R(\delta,f) + ||f||_\infty \cdot R(\delta,g)$.
    \end{enumerate}
\end{theorem}
\begin{proof}
    \begin{itemize}
        \item[(i)] For $\lambda \in \mathbb{R}$, we have 
        \begin{align*}
        R(\delta,\lambda f) &=\sum\limits_{S \in \Delta_{\delta}} R_{\lambda f}(S)\\
        &=\sum\limits_{S \in \Delta_{\delta}} \displaystyle \sup_{x,y \in S} |\lambda f(x) -\lambda f(y)|\\
        &= |\lambda| \sum\limits_{S \in \Delta_{\delta}} \displaystyle \sup_{x,y \in S} |f(x) - f(y)|\\
        &= |\lambda| \sum\limits_{S \in \Delta_{\delta}} R_{f}(S)= |\lambda| \cdot R(\delta, f).
        \end{align*}
        \item[(ii)] For the second one the proof is as follows:
        \begin{align*}
            R(\delta,f+g) &=\sum\limits_{S \in \Delta_{\delta}} R_{f+g}(S)\\
            &= \sum\limits_{S \in \Delta_{\delta}} \sup\limits_{x,y \in S} |(f+g)(x) - (f+g)(y)|\\
            & \le  \sum\limits_{S \in \Delta_{\delta}} \sup\limits_{x,y \in S} (|f(x) - f(y)|+|g(x)-g(y)|)\\
            &=  \sum\limits_{S \in \Delta_{\delta}} R_{f}(S) +  \sum\limits_{S \in \Delta_{\delta}} R_{g}(S)\\
            &= R(\delta,f) + R(\delta,g).
        \end{align*}
        \item[(iii)] The lines of proof for the last one is as follows
        \begin{align*}
            R(\delta,fg) &= \sum\limits_{S \in \Delta_{\delta}} R_{fg}(S)\\
            &= \sum\limits_{S \in \Delta_{\delta}} \sup\limits_{x,y \in S} |(fg)(x) - (fg)(y)|\\
            &=  \sum\limits_{S \in \Delta_{\delta}} \sup\limits_{x,y \in S} |f(x)g(x) - f(y)g(x) + f(y)g(x) - f(y)g(y)|\\
            &\le \sum\limits_{S \in \Delta_{\delta}}  \{ \sup\limits_{x,y \in S} |f(x) - f(y)||g(x)| +  \sup\limits_{x,y \in S} |f(y)||g(x) - g(y)|\}\\
            &\le  \sum\limits_{S \in \Delta_{\delta}}  \{ ||g||_\infty \sup\limits_{x,y \in S} |f(x) - f(y)| + ||f||_\infty \sup\limits_{x,y \in S} |g(x) - g(y)|\}\\
            &= ||g||_\infty \sum\limits_{S \in \Delta_{\delta}} R_{f}(S) + ||f||_\infty  \sum\limits_{S \in \Delta_{\delta}} R_{g}(S)\\
            &= ||g||_\infty \cdot R(\delta,f) + ||f||_\infty \cdot R(\delta,f).
        \end{align*}        
    \end{itemize}
\end{proof}
For $\beta \ge 1$, we define a function space
$$V_\beta [0,1] =\{ g \in \mathcal{C}[0,1] : ||g||_{\beta} < \infty \},$$
where $||g||_{\beta} = ||g||_\infty + \displaystyle \sup_{0 < \delta \le 1} \dfrac{\displaystyle R(\delta,g)}{\delta^{1-\beta}}$.
\begin{proposition}
    Let $f \in V_\beta [0,1]$, we define $||f||_\beta = ||f||_\infty +\displaystyle \sup_{0 < \delta \le 1} \dfrac{R(\delta , f)}{\delta^{1-\beta}}$. Then $||.||_\beta$ forms a norm on $V_{\beta}[0,1]$.
\end{proposition}
\begin{proof}
    \begin{enumerate}
        \item Let $f=0$. Then $||f||_\beta = ||0||_\beta = 0$.\\
        Conversely, let $||f||_\beta = 0.$ Then $||f||_\infty + \displaystyle \sup_{0 < \delta \le 1} \dfrac{R(\delta ,f)}{\delta^{1-\beta}} = 0$ \\
        $\implies ||f||_\infty = 0\ \text{and} \displaystyle \sup_{0 < \delta \le 1} \dfrac{R(\delta ,f)}{\delta^{1-\beta}} = 0 \implies f=0.$
        \item Let $\lambda (\neq 0) \in \mathbb{R}$.Then
        \begin{align*}
            ||\lambda f||_\beta &= ||\lambda f||_\infty + \sup_{0 < \delta \le 1} \dfrac{R(\delta ,\lambda f)}{\delta^{1-\beta}}\\
            &= |\lambda| \cdot ||f||_\infty + |\lambda| \cdot \sup_{0 < \delta \le 1} \dfrac{R(\delta ,f)}{\delta^{1-\beta}}\\
            &= |\lambda| \cdot ||f||_\beta.
        \end{align*}
        \item Let $f,g \in V_{\beta}(I)$. Then
        \begin{align*}
            ||f+g||_\beta &= ||f+g||_\infty + \sup_{0 < \delta \le 1} \dfrac{R(\delta ,f+g)}{\delta^{1-\beta}}\\
            &\le ||f||_\infty + ||g||_\infty + \sup_{0 < \delta \le 1} \dfrac{R(\delta ,f)}{\delta^{1-\beta}} + \sup_{0 < \delta \le 1} \dfrac{R(\delta ,g)}{\delta^{1-\beta}}\\
            &= ||f||_\beta + ||g||_\beta.
        \end{align*}
        \end{enumerate}
\end{proof}
\begin{lemma}
Let $\beta \in [1,2]$. Then $(V_\beta [0,1], ||.||_{\beta})$ is a Banach space.
\end{lemma}
\begin{proof}
    The lemma follows from Lemma 3.1. of \cite{16}.
\end{proof}
Now, we define non-stationary fractal functions on  the space $(V_\beta [0,1], ||.||_{\beta})$. For notational simplicity we denote the space $V_\beta [0,1]$ by $\mathcal{V}(I)$, where $I=[0,1]$.
\begin{theorem}
Let $f \in \mathcal{V}(I)$ and define
\[ \mathcal{V}_{f}(I) = \{ g \in \mathcal{V}(I) : g(x_0) = f(x_0), g(x_N) = f(x_N) \} .\]
Let $ b_m \in \mathcal{V}_{f}(I)$ be such that $\|b\|_\beta:=\sup\limits_{m\in \mathbb{N}}\|b_m\|_\beta<\infty$. Also, assume that the scaling functions $\alpha_{i,m}$ are constants such that $|\alpha|_\infty =\sup\limits_{m \in \mathbb{N}} \{ |\alpha_m|_{max} \} = \sup\limits_{m \in \mathbb{N}} \{ \max\limits_{i \in \mathbb{N}_N} |\alpha_{i,m}|\} < 1 $. Then the following hold.
\begin{enumerate}
\item The RB operator $T^{\alpha_m}$ defined in equation \eqref{eq-RB} is well defined on $\mathcal{V}_{f}(I)$.
\item In fact, $T^{\alpha_m} : \mathcal{V}_{f}(I) \longrightarrow \mathcal{V}_{f}(I) \subset \mathcal{V}(I)$ is a contraction map.
\item There exists a unique function $f_{b,V}^{\alpha} \in \mathcal{V}_{f}(I)$ such that the sequence\\
$\{ T^{\alpha_1}\ o\ T^{\alpha_2}\ o \dots o\ T^{\alpha_m} g \}$ converges to the map $f_{b,V}^{\alpha}$ for every $g \in \mathcal{V}_{f}(I)$.
\end{enumerate}
\end{theorem}
\begin{proof}
\begin{enumerate}
\item We have, 
\begin{align*}
    ||T^{\alpha_m}g||_{\beta} &= ||T^{\alpha_m} g||_\infty + \displaystyle \sup_{0 < \delta \le 1} \dfrac{R(\delta ,T^{\alpha_m}g)}{\delta^{1-\beta}}\\
    &= ||T^{\alpha_m} g||_\infty + \displaystyle \sup_{0 < \delta \le 1} \dfrac{\sum\limits_{S \in \Delta_{\delta}} R_{T^{\alpha_m}g}(S)}{\delta^{1-\beta}}
\end{align*}
Now,
\begin{align*}
&R_{T^{\alpha_m} g}(S)\\
&= \sup_{x,y \in S} |T^{\alpha_m} g(x) - T^{\alpha_m} g(y) |\\
&\le \sup_{x,y \in S} |f(x) - f(y)| + \max_{i \in \mathbb{N}_{N}} \sup_{x,y \in S_i} | \alpha_{i,m}. (g - b_{m})(Q_{i}(x)) - \alpha_{i,m}. (g - b_{m})(Q_{i}(y))|,\\
&\hspace{4cm} \text{where} \; S_i \; \text{is}\; \text{a} \; \text{subset}\; \text{of}\; I_i\\
&= \sup_{x,y \in S} |f(x) - f(y)|\\
&\hspace{1cm} + \max_{i \in \mathbb{N}_{N}} \Big(|\alpha_{i,m}|\Big) \sup_{x,y \in S_i} |\left[(g(Q_{i}(x))-g(Q_{i}(y)))- (b_{m}(Q_{i}(x)) - b_{m}(Q_{i}(y)))\right]|\\
&\le \sup_{x,y \in S} |f(x) - f(y)| + \max_{i \in \mathbb{N}_{N}} \Big(|\alpha_{i,m}|\Big) \sup_{\tilde{x},\tilde{y} \in S} (|g(\tilde{x})-g(\tilde{y})|+|b_m(\tilde{x})-b_m(\tilde{y})|)\\
&= R_f(S) + \max_{i \in \mathbb{N}_{N}} \Big(|\alpha_{i,m}|\Big) ( R_g(S) + R_{b_m}(S)).
\end{align*}
Taking sum over $\Delta_{\delta}$, we get
\begin{align*}
    &\displaystyle \sum_{S \in \Delta_{\delta}} R_{T^{\alpha_m} g}(S) \le \displaystyle \sum_{S \in \Delta_{\delta}} R_f(S) + \max_{i \in \mathbb{N}_{N}} \Big(|\alpha_{i,m}|\Big) \left(\displaystyle \sum_{S \in \Delta_{\delta}} R_g(S) +\displaystyle \sum_{S \in \Delta_{\delta}} R_{b_m}(S)\right)\\
    &\implies \displaystyle \sup_{0 < \delta \le 1} \dfrac{\displaystyle \sum_{S \in \Delta_{\delta}} R_{T^{\alpha_m} g}(S)}{\delta^{1-\beta}}\\
    &\le \displaystyle \sup_{0 < \delta \le 1} \dfrac{\displaystyle \sum_{S \in \Delta_{\delta}} R_f(S)}{\delta^{1-\beta}} + \max_{i \in \mathbb{N}_{N}} \Big(|\alpha_{i,m}|\Big) \left( \displaystyle \sup_{0 < \delta \le 1} \dfrac{\displaystyle \sum_{S \in \Delta_{\delta}} R_g(S)}{\delta^{1-\beta}} + \displaystyle \sup_{0 < \delta \le 1} \dfrac{\displaystyle \sum_{S \in \Delta_{\delta}} R_{b_m}(S)}{\delta^{1-\beta}}\right).
\end{align*}

Therefore, using the definition of the norm $||.||_\beta$, we get
$$||T^{\alpha_m}g||_{\beta}-||T^{\alpha_m} g||_\infty \le ||f||_{\beta}-||f||_\infty +\max_{i \in \mathbb{N}_{N}} \Big(|\alpha_{i,m}|\Big) (||g||_{\beta}-||g||_\infty + ||b_m||_{\beta}-||b_m||_\infty) $$
$$\implies ||T^{\alpha_m}g||_{\beta} \le ||f||_{\beta} + \max_{i \in \mathbb{N}_{N}} \Big(|\alpha_{i,m}|\Big) (||g||_{\beta} + ||b_m||_{\beta} ). $$

Since $f,g,b_m \in \mathcal{V}(I)$, the previous estimate ensures that $||T^{\alpha_m}g||_{\beta} < \infty $ and hence that $T^{\alpha_m}g \in \mathcal{V}(I)$.\\
Also $T^{\alpha_m}g(x_0) = f(x_0)$ and $T^{\alpha_m}g(x_N) = f(x_N)$, so that $T^{\alpha_m}g \in \mathcal{V}_{f}(I)$.
Therefore the RB operator is well defined on $\mathcal{V}_{f}(I)$.

\item Let $g_1, g_2 \in \mathcal{V}_{f}(I)$. For $ x \in I_i$,
\begin{equation*}
\begin{split}
|(T^{\alpha_m}g_1 - T^{\alpha_m}g_2)(x)| &= |\alpha_{i,m}||(g_1 - g_2)(Q_i(x))|\\
&\le \max_{i \in \mathbb{N}_{N}} \Big(|\alpha_{i,m}|\Big) ||g_1 - g_2||_\infty,
\end{split}
\end{equation*}
and hence
$$||(T^{\alpha_m}g_1 - T^{\alpha_m}g_2)||_\infty \le \max_{i \in \mathbb{N}_{N}} \Big(|\alpha_{i,m}|\Big) \cdot ||g_1 - g_2||_\infty .$$
Along lines similar to the estimation of $R_{T^{\alpha_m} g}(S)$, we get 
$$R_{T^{\alpha_m}g_1-T^{\alpha_m}g_2}(S) \le \max_{i \in \mathbb{N}_{N}} \Big(|\alpha_{i,m}|\Big) R_{g_1 - g_2} (S).$$
Therefore,
\begin{equation*}
\displaystyle \sup_{0 < \delta \le 1} \dfrac{\displaystyle \sum_{S \in \Delta_{\delta}} R_{T^{\alpha_m}g_1-T^{\alpha_m}g_2}(S)}{\delta^{1-\beta}} \le \max_{i \in \mathbb{N}_{N}} \Big(|\alpha_{i,m}|\Big) \displaystyle \sup_{0 < \delta \le 1} \dfrac{\displaystyle \sum_{S \in \Delta_{\delta}} R_{g_1 - g_2} (S)}{\delta^{1-\beta}}.
\end{equation*}
\begin{align*}
    \implies& ||T^{\alpha_m}g_1 - T^{\alpha_m}g_2||_{\beta}- ||(T^{\alpha_m}g_1 - T^{\alpha_m}g_2)||_\infty\\
    &\le \max_{i \in \mathbb{N}_{N}} \Big(|\alpha_{i,m}|\Big) \Big( ||g_1 - g_2||_{\beta} - ||g_1 - g_2||_\infty \Big).
\end{align*}
\begin{align*}
||T^{\alpha_m}g_1 - T^{\alpha_m}g_2||_{\beta} &\le \max_{i \in \mathbb{N}_{N}} \Big(|\alpha_{i,m}|\Big) (||g_1 - g_2||_{\beta} - ||g_1 - g_2||_\infty)\\
& \hspace{0.3 cm}+ ||(T^{\alpha_m}g_1 - T^{\alpha_m}g_2)||_\infty\\
&\le \max_{i \in \mathbb{N}_{N}} \Big(|\alpha_{i,m}|\Big) (||g_1 - g_2||_{\beta} - ||g_1 - g_2||_\infty)\\
&\hspace{0.3 cm}+ \max_{i \in \mathbb{N}_{N}} \Big(|\alpha_{i,m}|\Big) |g_1 - g_2||_\infty\\
&\le |\alpha|_\infty ||g_1 - g_2||_{\beta}.
\end{align*}
The assumption on the scaling sequence ensures that $T^{\alpha_m}$ is a contraction on $\mathcal{V}_{f}(I)$.
\item Let $g \in \mathcal{V}_{f}(I)$ be arbitrary. We show that $\{ ||T^{\alpha_m} g -g||_{\beta}\}$ is bounded. Adapting similar calculation as in the estimation of $ R_{T^{\alpha_m}g}(S)$, we get
$$R_{T^{\alpha_m}g - g}(S) \le R_f(S) + \max_{i \in \mathbb{N}_{N}} \Big(|\alpha_{i,m}|\Big)  R_{b_m}(S).$$
Thus,
$$\sup\limits_{0 < \delta \le 1} \dfrac{\sum\limits_{S \in \Delta_{\delta}} R_{T^{\alpha_m}g - g}(S)}{\delta^{1-\beta}} \le \sup\limits_{0 < \delta \le 1} \dfrac{ \sum\limits_{S \in \Delta_{\delta}} R_f(S)}{\delta^{1-\beta}} + \max_{i \in \mathbb{N}_{N}} \Big(|\alpha_{i,m}|\Big) \sup\limits_{0 < \delta \le 1} \dfrac{ \sum\limits_{S \in \Delta_{\delta}} R_{b_m}(S)}{\delta^{1-\beta}}.$$
$$\implies ||T^{\alpha_m}g -g||_{\beta}- ||(T^{\alpha_m}g -g)||_\infty \le ||f||_{\beta}-||f||_\infty + \max_{i \in \mathbb{N}_{N}} \Big(|\alpha_{i,m}|\Big) (||b_m||_{\beta}-||b_m||_\infty).$$
$\implies ||T^{\alpha_m}g -g||_{\beta} \le ||f||_{\beta} + \displaystyle \max_{i \in \mathbb{N}_{N}} \Big(|\alpha_{i,m}|\Big) ||b_m||_{\beta} \le ||f||_{\beta} + |\alpha|_\infty ||b||_{\beta}$.\\
 Applying Theorem \ref{mainthm}, $\exists$ a unique $f_{b,V}^{\alpha} \in \mathcal{V}_{f}(I)$ such that \\$f_{b,V}^{\alpha} = \displaystyle \lim_{ m \rightarrow \infty} T^{\alpha_1}\ o\ T^{\alpha_2}\ o \dots o\ T^{\alpha_m} g$ for any $g \in \mathcal{V}_{f}(I)$.
\end{enumerate}
\end{proof}
\begin{remark}
    In \cite{16}, Falconer and Fraser proved that for $\beta \in [1,2), C_\beta[0,1] = \bigcap\limits_{k \in \mathbb{N}} V_{\beta + \dfrac{1}{k}}[0,1].$ It is not known whether $C_\beta[0,1] $ is a complete normed space. A suitable norm on $C_\beta [0,1]$  may attract the researchers to define non-stationary fractal function on $C_\beta [0,1]$ and calculate its dimension. That is, for any $\beta \in [1,2)$, one can find a fractal function with dimension less than or equal to $\beta$ with respect to a suitable norm.

    We will construct a fractal function on the convex Lipschitz space and calculate its dimension in the next subsection.
\end{remark}
\subsection{Convex Lipschitz space}
\begin{definition}
Let $I = [a, b]$ and $\theta : \mathbb{R}^+ \longrightarrow \mathbb{R}^+$. A function $g$ is called convex Lipschitz of order $\theta$ on an interval $I$ provided there exists a constant $M$ such that $$|\Delta(u, v, \delta)|:= | g(u + \delta v) - (\delta g(u + v) + (1 - \delta)g(u))| \leq M \theta(v),$$
for $a \leq u < u + v \leq b$ and $0 \leq \delta \leq 1$.
\end{definition}
The set of all convex Lipschitz functions of order $\theta$ on $I$ is denoted by $\mathcal{V}^\theta(I)$. That is,
$$\mathcal{V}^\theta(I) = \{g : I \longrightarrow \mathbb{R} : g\ \text{is convex Lipschitz of order}\ \theta \}.$$
It is simple to verify that $\mathcal{V}^\theta$ is a vector space over the field $\mathbb{R}$, which we call the convex Lipschitz space of order $\theta$.
For $g \in \mathcal{V}^\theta(I)$, we define
$||g||_{\mathcal{V}^\theta} = ||g||_\infty + [g]^*$, where
$$[g]^* = \sup\limits_{a \le u < u+v \le b}
\dfrac{|\Delta(u, v, \delta)|}{\theta(v)} =\sup\limits_{a \le u < u+v \le b}
\dfrac{| g(u + \delta v) - (\delta g(u + v) + (1-\delta)g(u))|}{\theta(v)}.$$
It is simple to verify that $||.||_{\mathcal{V}^\theta}$ defines a norm on $\mathcal{V}^\theta(I)$.

\begin{theorem}
Let $f$ be a convex Lipschitz function of order $\theta$. Then
\begin{enumerate}
    \item $[\alpha f]^* = |\alpha| [f]^*$
    \item $[f \pm g]^* \le [f]^* + [g]^*$.
\end{enumerate}
\end{theorem}
\begin{proof}
It follows from the definition of $[f]^*$.
\end{proof}
\begin{proposition} \cite{18}
The convex Lipschitz space $\mathcal{V}^\theta(I)$ with respect to the norm $||.||_{\mathcal{V}^\theta}$ forms a complete metric space.
\end{proposition}

Our following proposition is collected from \cite{M}; it will help us to calculate the dimension of the non-stationary fractal function, which will occur in the upcoming theorem.

\begin{proposition}\label{prop}
Let $\theta : \mathbb{R}^+ \longrightarrow \mathbb{R}^+$ be a continuous map such that:
\begin{enumerate}
    \item For $s > 0,\ \theta(s) > 0 ;$
    \item $\limsup\limits_{s \rightarrow 0} [\dfrac{s}{\theta(s)}] < \infty $ and
    \item there exists $\gamma \ge 0$ such that $\lim\limits_{s \rightarrow 0} \frac{\theta(cs)}{\theta(s)} = c^\gamma$ for all $c > 0.$
\end{enumerate}
If $g \in \mathcal{C}[0,1] \cap \mathcal{V}^\theta(I)$, then 
\begin{enumerate}
\item for $\theta(s) = s^\epsilon$, $dim_H (Graph\ (g)) \le \overline{dim}_B (Graph\ (g)) \le 2-\epsilon.$

\item for $\theta(s) = -s\ln{s}$, $dim_H (Graph\ (g)) = dim_B (Graph\ (g)) = 1.$
\end{enumerate}
\end{proposition}

\begin{theorem}\label{convex}
Suppose that $f \in \mathcal{V}^\theta(I)$ and define
$$\mathcal{V}_{f}^{\theta}(I) = \{ g \in \mathcal{V}^\theta(I) : g(x_0) = f(x_0), g(x_N) = f(x_N) \}. $$
Suppose $ b_m \in \mathcal{V}_{f}^{\theta}(I)$ be such that $[b]^*: = \sup\limits_{m \in \mathbb{N}} [b_m]^* < \infty$ and the scaling functions $\alpha_{i,m}$ are constants such that
$S:= \max \{\max\limits_{i}|\alpha_{i,m}|,\; \max\limits_{i} |\alpha_{i,m}| \dfrac{\theta(Y)}{\theta(a_i Y)} \} < 1,$ where $Y=\frac{y}{a_i}$. Then the following hold.
\begin{enumerate}
\item The RB operator $T^{\alpha_m}$ defined in equation \eqref{eq-RB} is well defined on $\mathcal{V}_{f}^{\theta}(I)$.
\item In fact, $T^{\alpha_m} : \mathcal{V}_{f}^{\theta}(I) \longrightarrow \mathcal{V}_{f}^{\theta}(I) \subset \mathcal{V}^{\theta}(I)$ is a contraction map.
\item There exists a unique function $f_{b,\mathcal{V}^\theta}^{\alpha} \in \mathcal{V}_{f}^{\theta}(I)$ such that for every $g \in \mathcal{V}_{f}^{\theta}(I)$ the sequence $\{ T^{\alpha_1}\ o\ T^{\alpha_2}\ o \dots o\ T^{\alpha_m} g \}$ converges to the map $f_{b,\mathcal{V}^\theta}^{\alpha}$ .
\end{enumerate}
\end{theorem}

\begin{proof}
\begin{enumerate}
\item Clearly, $ \mathcal{V}_{f}^{\theta}(I) \subset \mathcal{V}^{\theta}(I)$ is a closed subset of $\mathcal{V}^{\theta}(I)$ and so that with respect to the metric induced by the norm $||.||_{\mathcal{V}^\theta}$,  $\mathcal{V}_{f}^{\theta}(I)$ is a complete metric space.
We have, $$(T^{\alpha_m}g)(x) = f(x) + \alpha_{i,m}. (g-b_{m})(Q_{i}(x)),$$
where $Q_{i}(x) = l_{i}^{-1}(x) = (x-e_i)/a_i $.\\
As $f,g,b_m \in \mathcal{V}_{f}^{\theta}(I), (T^{\alpha_m}g)(x) \in \mathcal{V}_{f}^{\theta}(I)$. 
Therefore the RB-operator is well defined on $\mathcal{V}_{f}^{\theta}(I)$.
\item For $g,h \in \mathcal{V}_{f}^{\theta}(I)$, \[||T^{\alpha_m}g - T^{\alpha_m}h||_{\mathcal{V}^\theta} = ||T^{\alpha_m}g - T^{\alpha_m}h||_\infty + [T^{\alpha_m}g - T^{\alpha_m}h]^*.\]
Let $Q_i(x) = X\ \text{and}\ y/a_i =Y$.
Now, 
\begin{align*}
    &[T^{\alpha_m}g - T^{\alpha_m}h]^*\\
    =&\sup\limits_{a \le x < x + y \le b} \left[\dfrac{| (T^{\alpha_m}g - T^{\alpha_m}h)(x + \delta y) - (\delta (T^{\alpha_m}g - T^{\alpha_m}h)(x + y) + (1-\delta)(T^{\alpha_m}g - T^{\alpha_m}h)(x))|}{\theta(y)}  \right]\\
   =& \max_{i} \sup\limits_{a \le a_i X + e_i < a_i (X + Y)+e_i \le b}
\Big[\dfrac{|\alpha_{i,m}|| (g-h)(X + \delta Y) - (\delta (g-h)(X + Y) + (1-\delta)(g-h)(X))|}{\theta(Y)} \times \\
& \hspace{7cm}\dfrac{\theta(Y)}{\theta(a_i Y)} \Big]\\
=& \max_{i} |\alpha_{i,m}| \times \dfrac{\theta(Y)}{\theta(a_i Y)} \times [g-h]^* .
\end{align*}

Also, \[||T^{\alpha_m}g - T^{\alpha_m}h||_\infty = |\alpha_{i,m}||(g-h)(Q_i(x))| \le \max_{i} |\alpha_{i,m}| \cdot ||g-h||_\infty.\]

Therefore, 
\begin{align*}
    ||T^{\alpha_m}g - T^{\alpha_m}h||_{\mathcal{V}^\theta} &\le \max_{i} |\alpha_{i,m}| \cdot ||g-h||_\infty + \max_{i} |\alpha_{i,m}| \times \dfrac{\theta(Y)}{\theta(a_i Y)} \times [g-h]^* \\
    &\le S \cdot ||g-h||_\infty + S \cdot [g-h]^*\\
    &= S ||g-h||_{\mathcal{V}^\theta}.
    \end{align*}

Since $S<1$, $T^{\alpha_m}$ is a contraction map for each $m \in \mathbb{N}$.

\item Let us take an arbitrary function $g \in \mathcal{V}^\theta(I)$. We have to check if the sequence $\{|| T^{\alpha_m}g - g||_{\mathcal{V}^\theta} \}$ is bounded.
$$||T^{\alpha_m}g - g||_{\mathcal{V}^\theta} = ||T^{\alpha_m}g - g||_\infty + [T^{\alpha_m}g - g]^*.$$
Now, 
\begin{align*}
    &[T^{\alpha_m}g - g]^*\\
    &=\sup\limits_{a \le x < x + y \le b} \left[\dfrac{| (T^{\alpha_m}g - g)(x + \delta y) - (\delta (T^{\alpha_m}g - g)(x + y) + (1-\delta)(T^{\alpha_m}g - g)(x))|}{\theta(y)}  \right] \\
    &\le \sup\limits_{a \le x < x + y \le b} \left[\dfrac{| (f - g)(x + \delta y) - (\delta (f - g)(x + y) + (1-\delta)(f - g)(x))|}{\theta(y)}  \right]\\ 
    & + \max_{i} \sup\limits_{a \le a_i X + e_i < a_i (X + Y)+e_i \le b}\\
&\left[\dfrac{|\alpha_{i,m}|| (g-b_m)(X + \delta Y) - (\delta (g-b_m)(X + Y) + (1-\delta)(g-b_m)(X))|}{\theta(Y)} \times \dfrac{\theta(Y)}{\theta(a_i Y)} \right] \\
&= [f-g]^* + \max_{i} |\alpha_{i,m}| \times \dfrac{\theta(Y)}{\theta(a_i Y)} \times [g-b_m]^* \\
&\le [f]^* + [g]^* + S \cdot ([g]^*+[b_m]^*) \\
&\le [f]^* + (1+S) \cdot [g]^* + S \cdot [b]^*.
\end{align*}
So the bound is independent of $m$. Applying Theorem \ref{mainthm}, there exists a unique $f_{b,\mathcal{V}^\theta}^{\alpha} \in \mathcal{V}_{f}^{\theta}(I)$ such that $f_{b,\mathcal{V}^\theta}^{\alpha} =  \lim\limits_{ m \rightarrow \infty} T^{\alpha_1}\ o\ T^{\alpha_2}\ o \dots o\ T^{\alpha_m} g$ for any $g \in \mathcal{V}_{f}^{\theta}(I)$.
\end{enumerate}
\end{proof}
\begin{remark}
    If we take $\alpha_{i,m} = \alpha_i$ for all $m \in \mathbb{N}$, then we get 
    $$f_{b,\mathcal{V}^\theta}^{\alpha} = \displaystyle \lim_{ m \rightarrow \infty} T^{\alpha}\ o\ T^{\alpha}\ o \dots o\ T^{\alpha} g = \displaystyle \lim_{ m \rightarrow \infty} (T^{\alpha})^m$$
    for any $g \in \mathcal{V}_{f}^{\theta}(I).$
    So that, $f_{b,\mathcal{V}^\theta}^{\alpha} \longrightarrow f_{\mathcal{V}^\theta}^{\alpha},\ \text{as}\ m \to \infty$, where $f_{\mathcal{V}^\theta}^{\alpha}$ is the stationary $\alpha$-fractal function on convex Lipschitz space $\mathcal{V}^\theta (I)$ that appeared in \cite{18}.
\end{remark}
\begin{theorem}
    Let $f, b_m (m \in \mathbb{N}) \in \mathcal{V}_{f}^{\theta}(I)$ and $\alpha_{i,m}$ are constants satisfying all the hypotheses of Theorem \ref{convex}.
    Also let $f$ be continuous on $[0,1]$. Then
    \begin{enumerate}
        \item for $\theta(t) = t^\epsilon$, $dim_H (Graph\ (f_{b,\mathcal{V}^\theta}^{\alpha})) \le \overline{dim}_B (Graph\ (f_{b,\mathcal{V}^\theta}^{\alpha})) \le 2-\epsilon.$
        
        \item for $\theta(t) = -t\ln{t}$, $dim_H (Graph\ (f_{b,\mathcal{V}^\theta}^{\alpha})) = dim_B (Graph\ (f_{b,\mathcal{V}^\theta}^{\alpha})) = 1.$
    \end{enumerate}
\end{theorem}

\begin{proof}
    For the given $S < 1$, we have from Theorem \ref{convex}, $f_{b,\mathcal{V}^\theta}^{\alpha} \in \mathcal{V}^\theta (I)$.\\
    We conclude the above results by applying Proposition \ref{prop}.
\end{proof}

\end{document}